\documentclass[12pt]{article} \textheight 9in \topmargin -.25in

\usepackage{amsmath, amsfonts, amsthm, amssymb, epsfig}

\begin{document}
\title{Pythagorean Triangles with Repeated Digits -- Different Bases}

\author{Habib Muzaffar\\
Department of Mathematics\\
International Islamic University Islamabad\\
P.O. Box 1244, Islamabad, Pakistan\\
\\
Konstantine Zelator\\
Department of Mathematics\\
College of Arts and Sciences\\
Mail Stop 942 \\
University of Toledo\\
Toledo, OH  43606-3390\\
USA}

\maketitle

\section{Introduction}

In 1998, in the winter issue of {\it Mathematics and Computer
Education}   (\cite{1}) Monte Zerger posed the following problem.
He had noticed or discovered the Pythagorean triple
$(216,630,666);\ (216)^2 + (630)^2 = (666)^2$. Note that $216=6^3$
and $666$ is the hypotenuse length of this Pythagorean triangle.
The question was, then whether there existed a digit $d$ (in the
decimal system) and a positive integer $k$ (other than the above)
such that $d^k$ is a leg length and $\underset{k\,{\rm
times}}{\left(\underbrace{d\ldots d}\right)}$ is the hypotenuse
length of a Pythagorean triangle.  The symbol or notation
$(\underset{k\,{\rm times}}{\underbrace{d \ldots d}})$ stands for
a natural number which in the base $10$ or decimal system has $k$
digits all of which are equal to $d$.  In other words
$(\underset{k\,{\rm times}}{\underbrace{d \ldots d}}) = d \cdot
10^{k-1} + d \cdot 10^{k-2} + \ldots + 10d + d$.

In 1999, F. Luca and Paul Bruckman (\cite{2}), answered the above
question in the negative.  They proved that the above Pythagorean
triple is the only one with this base $b=10$ property. In 2001, K.
Zelator took this question further and showed that there exists no
Pythagorean triangle one of whose legs having leg length $d^k$
while the other leg length being equal to $\underset{k\,{\rm
times}}(\underbrace{d \ldots d})$ (again with base $b=10$)
(\cite{3}). Note by the way, that any such triangle (of either
type) must be {\it non-primitive}.

The purpose of this work is to explore such questions in general,
when the base $b$ is no longer $10$.  In Section 2, we give
definitions and introduce some notation, while in Section 5 we
prove the five theorems of this paper.  We present a summary of
results in Section 3.  In Section 4, we state the very well-known
parametric formulas that generate the entire family of Pythagorean
triangles.  In Section 6, we give five families of Pythagorean
triangles with certain properties (similar to the triple
$(216,630,666)$ above; but with respect to bases $b$ other than
$10$); and in Section 7, we offer some corollaries to these
families.

Also, let us point out that in the proofs found in this work, we
only use elementary number theory.

\section{Notation and Definitions}

Let $b$ and $d$ be positive integers such that $b \geq 3$ and $2
\leq d \leq b-1$.

\vspace{.15in}

\noindent \fbox{\parbox{5.0in}{{\bf Notation:}  By $d_{k,b}$,
where $k$ is a positive integer, $k \geq 2$, we will mean the
positive integer which in the base $b$ system has $k$ digits all
equal to $d$.  In other words,
$$\begin{array}{rcl}{d_{k,b}}& =& d\cdot \left(b^{k-1} +
\ldots + b+1\right) = \dfrac{d\cdot (b^k-1)}{b-1}\\
& = & d\cdot b^{k-1} + d\cdot b^{k-2}+\ldots + d \cdot b +
d\end{array}$$
$d_{k,b}= (\underset{k\,{\rm
times}}{\underbrace{d\ldots d}})_b.$ Also, we denote a Type 1
triangle (see definition below) by $\mathbf{T_1(k,b,d)}$; and a
Type 2 triangle by $\mathbf{T_2(k,b,d)}$.  Note that for $k \geq
2,\ d_{k,b} > db^{k-1} > d \cdot d^{k-1}=d^k$.}}

\vspace{.15in}

\noindent \fbox{\parbox{5.0in}{{Definition 2:} \begin{enumerate}
\item[(i)] A Pythagorean triangle is called a {\it Type 1 triangle
with base $b$ repeated digits},
 if there exist positive integers $k,b,d$ such that $b\geq 3,\ 2 \leq
 d \leq b-1,\ k\geq 2$; and with one
of its two legs having length $d^k$, while the hypotenuse having
length $d_{k,b}$.  We denote such a triangle by $T_1(k,b,d)$.

\item[(ii)]  A Pythagorean triangle is called a {\it  Type 2
triangle with base $b$ repeated digits}, if there exist positive
integers $k,b,d$ such that $b \geq 3,\ 2 \leq d \leq b-1,\ k\geq
2$; and with one of its legs having length $d^k$, while the other
length being equal to $d_{k,b}$.  We denote such a triangle by
$T_2(k,b,d)$.
 \end{enumerate}}}

\vspace{.15in}

 \noindent {\bf Remarks}

 \begin{enumerate}
 \item[1.]  Note that the inequalities in the above definition are
 justified by inspection, by the fact that no Pythagorean triangle
 can have a side whose length is equal to $1$.  Consequently,
 $d\geq 2$ and thus $b$ must be at least $3$ in value.  Also, as
 trivially, one can see that $k$ must be at least $2$ in value.
 (No Pythagorean triangle can be isosceles.)

 \item[2.] Since two side lengths completely determine a
 Pythagorean triangle, it follows that both notations $T_1(k,b,d)$
 and $T_2(k,b,d)$ are unambiguous.  In other words, $T_1(k,b,d)$
 or $T_2(k,b,d)$ can only represent one Pythagorean triangle.
 \end{enumerate}

 \section{Summary of Results}

 In Theorem 1, we prove that if $b=4$, there exist no Type 2
 Pythagorean triangles.  In other words, there exists no
 Pythagorean triangle which is a $T_2(k,4,d)$, for some $k$ and
 $d$.  In Theorem 2, it is shown that the only Pythagorean
 triangle which is a $T_1(k,4,d)$ is the triangle $T_1(2,4,3)$,
 which has side  lengths $9,12,$ and $15$.  In Theorem 3, we prove
 that there exists no Pythagorean triangle which is a
 $T_2(k,3,d)$.  Likewise, for triangles $T_1(k,3,d)$.  No such
 triangle exists, according to Theorem 4.  Finally, Theorem 5 says
 that there exist no Pythagorean triangles of the form
 $T_2(2,b,d)$ with $2\leq d \leq 4$, and for any value of $b$
 (remember that always, $2 \leq d \leq b-1$).  Of the five
 families presented in Section 5, three are families of triangles
 which are $T_1(k,b,d)$.  The other two families consist of Type 2
 triangles or $T_2(k,b,d)$.

\section{Pythagorean triples-parametric formulas}

If $(a,b,c)$ is a Pythagorean triple with $c$ being the hypotenuse
length then (without loss of generality; $a$ and $b$ can be
switched)

\vspace{.15in}

\noindent \fbox{\parbox{5.0in}{$a = \delta(m^2-n^2), \ b =
\delta(2mn),\ c= \delta(m^2+n^2)$, where $\delta,m,n$ are positive
integers such that $m > n \geq 1$,  $(m,n)=1$ (i.e., $m$ and $n$
are relatively prime) and $m+n \equiv 1({\rm mod\,} 2)$ (i.e., $m$
and $n$ have different parities; one is even, the other odd). }}
\hfill (1)

\vspace{.15in}

The above formulas are very well known, and they generate the
entire family of Pythagorean triangles.  They can be found in
almost every number theory book; certainly in any undergraduate
number theory textbook.

\section{The Five Theorems and Their Proofs}

\noindent \fbox{\parbox{5.0in}{{\bf Theorem 1:}

There exists no Pythagorean triangle which is a $T_2(k,4,d)$.}}

\vspace{.15in}

\noindent {\bf Proof:}  If such a triangle existed, we would have
$b=4,k\geq 2$, and $2 \leq d \leq 3$.  According to (1),

\vspace{.15in}

$\left.\left.\begin{array}{lll}{\rm either} & d_{k,4} = \delta
(m^2-n^2), &
d^k = \delta(2mn) \\
{\rm or} & d^k = \delta (m^2-n^2), & d_{k,4} = \delta (2mn)
\end{array}\right\}  \begin{array}{l}(2a) \\ (2b)\end{array} \right\}$ \hfill
(2)

\vspace{.15in}

\noindent We distinguish between cases, according to whether $d=2$
or $d=3$.

\vspace{.15in}

\noindent {\it Case 1: d=2}

Since $m$ and $n$ have different parities, it follows that $mn
\equiv 0({\rm mod\,}2)$ and thus $\delta \cdot (2mn) \equiv 0
({\rm mod\,}4)$.  Also,

$$
d_{k,4} = 2_{k,4} = \dfrac{2\cdot (4^k-1)}{4-1} =
\dfrac{2(4^k-1)}{3} \equiv 2({\rm mod\,} 4),$$

\noindent which shows that possibility (2b) is ruled out.  Thus,
we consider (2a) with $d=2$:

\vspace{.15in}

\hspace{1.0in} $\dfrac{2(4^k-1)}{3} = \delta (m^2-n^2),\ \ 2^k=
\delta(2mn)$ \hfill (2c)

\vspace{.15in}

\noindent The second equation (2c) implies, by virtue of $m > n$,\
$(m,n)=1,$ and $m+n \equiv 1({\rm mod\,} 2)$, that

\hspace{1.0in} $\left.\begin{array}{l} m=2^u,\ n=1,\ \delta =
2^v;\ {\rm
for\ some}\\
{\rm integers}  \ \ u \ {\rm and}\ v\ {\rm such \ that\ } v \geq
0, u \geq 1 \\
{\rm and\ with}\ u+v+1=k \end{array}\right\}$ \hfill (2d)

\vspace{.15in}  By the first equation in (2c) and (2d) we obtain,

\vspace{.15in}

\hspace{1.0in} $2 \cdot (4^k-1) = 3\cdot 2^v(2^{2u}-1) $ \hfill
(2e)

\vspace{.15in}

\noindent which easily implies $ v=1$ (consider the power of $2$
on both sides).  Thus, from $k=u+v+1 \Rightarrow k=u+2$; and by
(2e),

$$
2^{2u+4} + 2 = 3 \cdot 2^{2u} \Leftrightarrow 2^{2u+3} + 1 = 3
\cdot 2^{2u-1},
$$

\noindent which is impossible modulo 2 since $2u-1 \geq 1$ (by
(2d)).

\vspace{.15in}

\noindent {\it Case 2: d=3}

We have $d_{k,4} = 3_{k,4} = \dfrac{3(4^k-1)}{3} = 4^k-1$.
Obviously, since $4^k-1$ is odd, it cannot equal $\delta(2mn)$.
Thus again, as in the previous case, this leads us to (2a) in (2);
$$
4^k-1=\delta(m^2-n^2),\ \ 3^k = \delta(2mn),
$$

\noindent an impossibility again since $3^k$ is odd, while
$\delta(2mn)$ is even.  The proof is complete. \hfill $\Box$

\vspace{.15in}

\noindent \fbox{\parbox{5.0in}{{\bf Theorem 2:}

The only Type 1 triangle which is a $T_1(k,4,d)$, is the triangle
$T_1(2,4,3)$, which has side lengths $9,12,$ and $15$.}}

\vspace{.15in}

\noindent {\bf Proof:}  Let $T_1(k,4,d)$ be such a triangle; it's
hypotenuse length being $d_{k,4}$.  We must have

\vspace{.15in}

 $\begin{array}{l} {\rm either} \\ {\rm
or}\end{array}\left. \left\{ \begin{array}{ll} d^k =
\delta(m^2-n^2), & d_{k,4} = \delta (m^2+n^2)) \\ d^k =
\delta(2mn), & d_{k,4} = \delta (m^2+n^2) \end{array} \right\}
\begin{array}{l} (3a) \\ (3b)\end{array} \right\}$ \hfill (3)

\vspace{.15in}

\noindent Since $b=4$, we must have $d=2$ or $3$.  We distinguish
between two cases.

\vspace{.15in}

\noindent {\it Case 1: d=3}

Obviously, possibility (3b) cannot hold true since $3^k \not\equiv
0 ({\rm mod\,} 2)$.  Thus, we need only consider (3a).  Since
$d_{k,4} = 3_{k,4} = \dfrac{3(4^k-1)}{3} = 4^k -1$; equations (3a)
imply

\vspace{.15in}

$\left. \begin{array}{ll} & \delta (m^2-n^2) = \delta (m-n) (m+n)
= 3^k \\
{\rm and} & \delta (m^2+n^2) = 4^k-1 \end{array}\right\}$ \hfill
(3c)

\vspace{.15in}

Since $(m,n)=1$ and $m,n$ have different parities, we conclude
that $(m-n,m+n) = 1$ and also $1 \leq m-n < m+n$.  This, then
combined with the first equation in (3c), implies

\setcounter{equation}{3}

\vspace{.15in}

\begin{equation}
\left. \begin{array}{l} \delta = 3^v, m-n=1, m+n=3^w, w+v = k \\
{\rm for \ integers}\ v,w\ {\rm with}\ v \geq 0\ {\rm and}\ w \geq
1
\end{array}\right\} \label{E4} \end{equation}

\noindent From (\ref{E4}) we obtain $m = \dfrac{3^w+1}{2}$ and $n
= \dfrac{3^w-1}{2}$ and thus, the second equation in (3c) gives

\begin{equation}
2(4^k-1) = 3^v \cdot \left[ 3^{2w} + 1\right] \label{E5}
\end{equation}

By virtue of the fact that  $4^3=64$; $4^3-1 \equiv 0 ({\rm mod\,}
9)$, the following three statements can be easily verified.

\vspace{.15in}

\noindent If $k \equiv 0({\rm mod\,} 3) \Rightarrow 4^k-1 \equiv 0
({\rm mod\,} 9)$

\noindent If $k \equiv 1({\rm mod\,} 3) \Rightarrow 4^k-1 \equiv
3({\rm mod\,} 9)$

\noindent If $k \equiv 2({\rm mod\,} 3) \Rightarrow 4^k-1 \equiv
6({\rm mod\,} 9)$

\vspace{.15in}

\noindent This then shows, by (\ref{E5}), that if $v \geq 2$, $k$
must be a multiple of $3$.

Accordingly, we distinguish between two subcases:  $v \geq 2$
being one subcase, while $v < 2$ (i.e., $v = 0$ or $1$) the other.

\vspace{.15in}

\noindent {\it Subcase 1a: $v \geq 2$}

In this subcase, we must have $k \equiv 0({\rm mod\,} 3)$ (see
above).  Since $k \geq 2$, if we consider (5) modulo (8), we see
that on account of

\vspace{.15in}

$3^{2w} \equiv 9^w\equiv 1({\rm mod\,} 8)$; we have

$$\begin{array}{rcl} 2(0-1) & \equiv & 3^v (1+1) ({\rm mod\,}8);\\
-2 & \equiv &  3^v \cdot 2 ({\rm mod\,}8);\\
-1 & \equiv & 3^v ({\rm mod\,} 4) \Rightarrow v \equiv 1 ({\rm
mod\,} 2)
\end{array}
$$

\noindent (Clearly, if $v$ were even, $3^v$ would be congruent to
$1({\rm mod\,} 4)$.  Thus, $v$ must be an odd integer.  The next
observation shows that $w$ must also be odd.  To see why, observe
that if $w$ were even; then $2w\equiv 0({\rm mod\,}4) \Rightarrow
3^{2w} \equiv 1({\rm mod\,}16)$, since $3^4 \equiv 1({\rm mod\,}
16)$. But then, if we consider (5) modulo 16 we see that

$$
2(0-1)=3^v(1+1)({\rm mod\,}16) \Rightarrow  14 \equiv 2 \cdot
3^v({\rm mod\,}16),
$$

\noindent which is impossible because $2 \cdot 3^v \equiv 6({\rm
mod\,} 16)$, in view of $v \equiv 1({\rm mod\,} 2)$.  Indeed, to
make this a bit more clear, put $v = 4t + 1$ or alternatively $v =
4t +3$ for some integer $t$.  If $v = 4t+1 \Rightarrow 2 \cdot 3^v
= 2 \cdot 3^{4t+1} \equiv 2 \cdot 3^{4t} \cdot 3({\rm mod\,} 16)$
so that $2 \cdot 3^v \equiv 2 \cdot 1 \cdot 3 \equiv 6({\rm
mod\,}16)$. If, on the other hand, $v=4t+3$, we have $2\cdot 3^v =
2 \cdot 3^{4t+3} \equiv 2 \cdot 3^{4t} \cdot 3^3 \equiv 2\cdot
1\cdot 11 \equiv 22 \equiv 6({\rm mod\,}16)$.  Therefore, both $v$
and $w$ must be odd; and so, by $k = v + w$ in (4), it follows
that $k$ must be even.  Thus, since $k \equiv 0({\rm mod\,}2)$ and
$k \equiv 0({\rm mod\,} 3)$ (see beginning of this subcase), it
follows that $k \equiv 0 ({\rm mod\,} 6)$.  Next, we apply
Fermat's Little Theorem for the prime $p=7$: $4^k -1 \equiv ({\rm
mod\,}7)$ and thus, by (\ref{E5}) we see that
$$3^{2w} +1 \equiv 0({\rm mod\,} 7)$$

\noindent which is impossible by virtue of the fact that $w$ is
odd.  Indeed, $w \equiv 1,3,$ or $5({\rm mod\,}6) \Rightarrow 2w
\equiv 2,0,$ or $4({\rm mod\,} 6)4$ and so

$$
\begin{array}{rcl}
3^{2w}+1 & \equiv & 3^2+1, 3^0+1\ {\rm or}\ 3^4+1\\
& \equiv & 9+1,\ 1+1 \ {\rm or}\ 4+1;\\
& \equiv & 3, 2,\ {\rm or}\ 5({\rm mod\,}7).\end{array}
$$

This concludes the proof of Subcase 1a. \hfill $\Box$

\vspace{.15in}

\noindent {\it Subcase 1b: $v < 2; v=0$ or $1$}

If $v=0$, then from (4) we have $w=k$ and hence from (5),

$$
2\cdot (4^k-1) = 3^{2k} + 1 \Rightarrow 2 \cdot 4^k = 3 \cdot
(3^{2k-1}+1),$$

\noindent which is impossible modulo 8; since $2 \cdot 4^k \equiv
0({\rm mod\,}8)$ (in view of $k \geq 2$), and $3^{2k-1} +1 \equiv
4({\rm mod\,} 8)$; and so $3 (3^{2k-1}+1) \equiv 3 \cdot 4 \equiv
4 ({\rm mod\,} 8)$  as well.  Note that when the exponent is odd,
say $2\rho +1$ (like in the case of $2k-1$); $3^{2\rho
+1}+1=3^{2\rho} \cdot 3+1 \equiv 1 \cdot 3+1 \equiv 4({\rm
mod\,}8)$. If $v=1$, equation (\ref{E4}) gives $k=w+1$ and by
(\ref{E5})

$$
2(4^k-1)=3 \cdot \left[ 3^{2(k-1)} + 1 \right] \Leftrightarrow
2^{2k+1} = 3^{2k-1} + 5.
$$

\noindent We claim that $k$ must equal $2$; for if to the contrary
$k\geq 3$, the last equation implies $2 \cdot 2^{2k} =
\dfrac{3^{2k}}{3} + 5; \ 6 = \left( \dfrac{9}{4}\right)^k +
\dfrac{15}{4^k} \Rightarrow \left(\dfrac{9}{4}\right)^k < 6$ which
is impossible since for $k \geq 3;\ \left(\dfrac{9}{4}\right)^k
\geq \left( \dfrac{9}{4}\right)^3 > 6$.  Thus, $k=2$, and so from
(\ref{E4}) we obtain $w=1$.

Altogether, $v=1=w$ and $k=2$.  By (\ref{E4}) $\Rightarrow \delta
= 3, \ m=2, \ n=1$, we obtain the Pythagorean triangle whose side
lengths are $9,12,$ and $15$; and since $\delta = 3$,  this is the
triangle $T_1(2,4,3)$.  This concludes the proof of Subcase 1b.
\hfill $\Box$

\vspace{.15in}

\noindent {\it Case 2: d=2}

Going back to (3), we easily see that possibility (3a) cannot hold
true, since the first equation in (3a) would imply that both
$\delta$ and $m^2-n^2$ are powers of $2$.  But $m^2-n^2$ is an odd
integer (since $m$ and $n$ have different parities) and $m^2-n^2
>1$.

Now, consider (3b).  We have,

\begin{equation}
2^k=\delta (2mn),\ \dfrac{2(4^k-1)}{3} = \delta (m^2+n^2)
\label{E6}
\end{equation}

\noindent and since $m > n \geq 1$, $(m,n)=1$ and $m+n \equiv
1({\rm mod\,} 2)$, the first equation in (\ref{E6}) implies

\begin{equation}
\left.\begin{array}{l} \delta = 2^u, \ m = 2^t,\ n = 1,\ k = u + t
+ 1\\
{\rm for \ integers}\ u, t,\ {\rm with}\ u \geq 0\ {\rm and}\ t
\geq 1 \end{array}\right\} \label{E7}
\end{equation}

Combining (\ref{E7}) with the second equation in (\ref{E6}) yields

$$
4^k-1 = 3 \cdot 2^{u-1} \cdot (2^{2t}+1). $$

\noindent The possibility $u = 0$ is impossible since $3 \cdot
(2^{2t} +1)$ is odd.

It follows that we must have $u=1$; and by (\ref{E7}), $t = k-2$.
Thus, the last equation above gives, after some algebra,
$4^{k-1}-1 = 3 \cdot 2^{2k-6}$. Recall that $k \geq 2$. Clearly,
the last equation requires $2k-6=0$ since its left-hand side is an
odd integer. Thus, $k=3$, which in turn implies $4^2-1=3;\ 15=3$,
a contradiction.  This concludes the proof of Theorem 2. \hfill
$\Box$

\vspace{.15in}

\noindent \fbox{\parbox{5.0in}{{\bf Theorem 3:}

Let $b=3$.  There exists no Pythagorean Type 2  triangle with base
3 repeated digits.  In other words, there exists no triangle which
is a $T_2(k,3,d)$.}}

\vspace{.15in}

\noindent {\bf Proof:}  First observe that since $2 \leq d \leq
b-1,\ 2\leq d \leq 3-1 \Rightarrow d=2$; and so, $d^k =2^k$ and
$d_{k,3} = 2_{k,3} = \dfrac{2(3^k-1)}{2} = 3^k-1$.  If such a
triangle exists, one leg will have length $2^k$, the other
$3^k-1$.  Thus, there are two possibilities.

\vspace{.15in}

\begin{equation}\left.\left.\begin{array}{lll} {\rm Either} & 3^k-1 = \delta
(2mn), &
2^k = \delta(m^2-n^2) \\
{\rm or} & 3^k-1 = \delta(m^2-n^2), & 2^k = \delta(2mn)\end{array}
\right\} \begin{array}{l} (8a)\\ (8b)
\end{array}\right\}
\end{equation}

\noindent {\it Case 1:} Assume possibility (8a) to hold.

From the second equation in (8a), it follows that both positive
integers $\delta$ and $m^2-n^2$ are powers of $2$.  However,
$m^2-n^2 = (m-n)(m+n) \geq 3$, on account of $m > n \geq 1$.  In
fact, $m^2-n^2 \geq (n+1)^2 - n^2 = 2n+1 \geq 3$.  Since $m,n$
have different parities, $m^2-n^2$ is an odd integer greater than
or equal to $3$.  Thus, it cannot equal to a power of $2$, which
renders the second equation in (8a) contradictory.

\vspace{.15in}

\noindent {\it Case 2:}  Assume possibility (8b)

The second equation implies that each of the positive integers
$\delta, m, n$, must be a power of $2$; and since $m > n \geq 1$
and $m+n\equiv 1({\rm mod\,}2)$, it follows that

\setcounter{equation}{8}

\begin{equation}
\left. \begin{array}{ll} & n = 1,\ m=2^v, \ \delta = 2^p \\
{\rm for \ integers} & v\ {\rm and}\ p\ {\rm such\ that}\ v \geq
1, p \geq 0 \end{array}\right\} \label{E9}
\end{equation}

Combining (\ref{E9})  with the first equation in (8b) yields,

\begin{equation}
3^k -1 = 2^p \cdot (2^{2v}-1) \label{E10}
\end{equation}

Consider equation (\ref{E10}) modulo 3.  We have,

$$2^{2v} -1 = 4^v -1 \equiv 1 - 1 \equiv 0({\rm mod\,}3)\Rightarrow
2^p \cdot (2^{2v}-1)\equiv 0({\rm mod\,}3)$$

\noindent while

$$
3^k -1 \equiv 0 -1 \equiv -1 \equiv 2 ({\rm mod\,}3).$$

We have a contradiction.  The proof is complete. \hfill $\Box$

\vspace{.15in}

\noindent \fbox{\parbox{5.0in}{ {\bf Theorem 4:} Let $b =3$. There
exists no Type 1 triangle with base 3 repeated digits.  In other
words, there exists no triangle which is a $T_1(k,3,d)$}}

\vspace{.15in}

\noindent {\bf Proof:}

Again, as in the previous proof of Theorem 3, observe that in view
of $2 \leq d \leq b-1$, we have $d=2$.  If a triangle $T_2(k,3,2)$
exists, it would be a Pythagorean triangle with the hypotenuse
having length \linebreak $d_{k,3} = \left(\dfrac{3^k-1}{2}\right)
\cdot 2 = 3^k-1$, and with one of the two legs having length
$2^k$. Which means that,

\begin{equation}
\left.\left. \begin{array}{ll} {\rm Either} & 2^k =
\delta(m^2-n^2),
3^k-1 = \delta (m^2+n^2) \\
{\rm or} & 2^k = \delta (2mn), 3^k-1 = \delta (m^2 + n^2)
\end{array}\right\}\begin{array}{l}(11a)\\ (11b)
\end{array}\right\} \label{E11}\end{equation}

The first possibility (11a) is ruled out at once by the first
equation in (11a), since $m^2-n^2$ is an odd integer and $m^2 -
n^2 \geq 3$; we already saw this in the proof of Theorem 3.

Next, let us consider (11b), the other possibility.  As in the
previous proof, in view of $m > n \geq 1$, we easily infer from
the first equation that

\begin{equation} \left. \begin{array}{ll} & \delta = 2^p, m=2^v,
n=1\\
{\rm for\ some\ integers} & p,v \ {\rm with}\ p \geq 0, v\geq 1\\
{\rm and}&  p + v + 1 = k.
\end{array}\right\} \label{E12}
\end{equation}

Recall that always $k \geq 2$.  From (\ref{E12}) and the second
equation in (11b), we obtain

\begin{equation} 3^{p+v+1} -1 = 2^p\cdot (2^{2v}+1) \label{E13}
\end{equation}

Consider (\ref{E13}) modulo $3$.  Since $2^{2v} +1 = 4^v +1 \equiv
1+1 \equiv 2({\rm mod\,}3)$, (\ref{E13}) implies $-1 \equiv
2^{p+1} ({\rm mod\,}3) \Leftrightarrow 2 \equiv 2^{p+1} ({\rm
mod\,} 3) \Leftrightarrow$ (since $(2,3)=1$) $1 \equiv 2^p({\rm
mod\,}3) \Leftrightarrow p\equiv 0 ({\rm mod\,}2)$.  $p$ must be
an even integer. Moreover, since the left-hand side of (\ref{E13})
is an even integer, while $2^{2v}+1$ is odd; $2^p$ must be even;
which means $p\geq 1$.  But $p$ is even, so we must have $p \geq
2$.

Since $p \geq 2$, the right-hand side of (\ref{E13}) must be a
multiple of $4$;

\vspace{.15in}

\hspace{1.0in} $3^{p+v+1} - 1 \equiv 0({\rm mod\,}4)$ \hfill (13a)

\vspace{.15in}

However, if $\ell$ is a positive integer, then $3^{\ell} -1 \equiv
0({\rm mod\,}4)$ if $\ell$ is even; while $3^{\ell}-1 \equiv
2({\rm mod\,} 4)$ if $\ell $ is odd, as it can be easily verified.
Thus, (13a) implies that the exponent $p+v+1$ must be an even
integer; and in view of $p \equiv 0({\rm mod\,}2)$, we see that
$v$ must be odd:

\vspace{.15in}

\hspace{1.0in} $ p\equiv 0({\rm mod\,}2),\ \ v \equiv 1({\rm
mod\,}2)$ \hfill (13b)

\vspace{.15in}

Accordingly by (13b), $p \equiv 0,2,$ or $4({\rm mod\,} 6)$, while
$v \equiv 1,3,$ or $5({\rm mod\,} 6)$

\noindent We will use this, by considering equation (13) modulo
$7$.  To facilitate this end, observe that, if $r$ is a positive
integer and $r \equiv i({\rm mod\,} 6)$, with $0\leq i \leq 5$,
then $2^r \equiv 2^i ({\rm mod\,}7)$ and $3^r \equiv 3^i ({\rm
mod\,}7)$. This observation leads to the following table:

\vspace{.15in}

$\begin{array}{|c|c|c|c|}\hline &&&\\
 {\rm Value\ of} \ p & {\rm
Value \
of}\ v & {\rm Value\ of}\ 2^p\cdot(2^{2v}+1)& {\rm Value\ of}\ 3^{p+v+1}-1\\

{\rm modulo}\ 6 & {\rm mod\,}\ 6 & {\rm mod\,}\ 7 & {\rm mod\,} \ 7 \\
\hline &&&\\
 0 & 1 & 5 & 1 \\ \hline
 &&&\\
 0 & 3 & 2 & 3 \\ \hline
&&&\\
 0 & 5 & 3 & 0\\ \hline
 &&& \\
 2 & 1 & 6 & 3 \\ \hline
 &&&\\
 2 & 3 & 1 & 0 \\ \hline
&&& \\
2 & 5 & 5 & 1 \\ \hline
&&&\\
 4 & 1 & 3 & 0 \\ \hline
 &&&\\
 4 & 3 & 4 & 1 \\ \hline
 &&&\\
 4 & 5 & 6 & 3 \\ \hline
\end{array}$

\vspace{.15in}

The results on the last two columns (columns 3 and 4) clearly
render (13) impossible modulo $7$;  a contradiction.  The proof is
complete.  \hfill $\Box$

\vspace{.15in}

\noindent\fbox{\parbox{5.0in}{{\bf Theorem 5:}

Let $k=2,\ 2 \leq d \leq 4$ (i.e., $d=2,3,$ or $4$).  There exists
no Type 2 triangle with base $b$ repeated digits.  In other words,
there exists no triangle which is a $T_2(2,b,d)$.  That is no
triangle which is a $T_2(2,b,2)$, a $T_2(2,b,3)$, or a
$T_2(2,b,4)$.}}

\vspace{.15in}

\noindent {\bf Note:}  Due to $2 \leq d \leq b-1$, when $d=2$, we
must have $b\geq 3$; when $d=3,\ b \geq 4$; while for $d=4,\ b
\geq 5$.

\vspace{.15in}

\noindent {\bf Proof:}  We give a simple proof, without making use
of parametric formulas (1).  If a triangle $T_2(2,b,d)$ exists,
with $2 \leq d \leq 4$, then it must have one leg length equal to
$d_{2,b}= d(b+1)$, while the other leg length is equal to $d^2$.
Thus, $\left[d(b+1)\right]^2 + (d^2)^2 = m^2$, for some positive
integer $m$; and so,

$$d^2 \cdot \left[ (b+1)^2 +d^2\right] = m^2.
$$

\noindent In the last equation, $d^2$ is a divisor of $m^2$, so
$d$ must be a divisor of $m$.  Put $m = d\cdot c$, for some
positive integer $c$, in order to obtain

\setcounter{equation}{13}

\begin{equation}
d^2+(b+1)^2 = c^2 \Leftrightarrow \left[c+(b+1)\right] \cdot
\left[ c-(b+1)\right] = d^2 \label{E14}
\end{equation}

Next, we use the conditions $2 \leq d \leq 4$ and $2 \leq d \leq
b-1$.  If $d =2$, (\ref{E14}) $\Rightarrow$ (since $c+b+1 >
c-(b-1)$) $c+b+1 = 4$ and $c-(b+1)=1$ which, in turn, gives
$c=\dfrac{5}{2}$, a contradiction since $c$ is an integer.  If
$d=3$, (\ref{E14}) $\Rightarrow c+b+1=9$ and $c-(b+1)=1$, which
gives $b=3$, contradicting $d \leq b-1$.

If $d=4$, we have either $c+b+1=8$ and $c-(b+1)=2$ or $c+b+1=16$
and $c-(b+1)=1$.  In the first case we obtain $b=2$, a
contradiction since $b \geq 5$.  In the second case, $c =
\dfrac{17}{2}$ a contradiction once more.  The proof is complete.
\hfill $\Box$

\section{Five Families with $k =2$}

A basic principle used for the construction of all five families
below, is the identity $(r^2-q^2)^2 + (2rq)^2 = (r^2+q^2)^2$.

\vspace{.15in}

\noindent \fbox{\parbox{5.0in}{{\it Basic Principle:} $(r^2-q^2)^2
+ (2rq)^2 = (r^2+q^2)^2$, for any positive integer $r$ and $q$.}}

\vspace{.15in}

\begin{enumerate}
\item[A.] {\it Two families of Type 2 triangles}

In both families below, we take $k=2$; and so $d_{k,b} = d_{2,b} =
d\cdot b + d = d(b+1)$.  To see how the first family comes about,
let $\ell,q$ be positive integers such that $\ell^2 \leq 2q^2-2$,
and let $r = q+\ell$.  Take $b=2rq-1,\ d=r^2-q^2$. Note that
$d=r^2-q^2 \geq (q+1)^2 -q^2=2q+1 \geq 3$ and also that, $b-1-d =
2rq -1 -1 -(r^2-q^2) = 2(q+\ell)q-2-(q+\ell)^2 + q^2 =
2q^2-\ell^2-2 \geq 0$.  Thus $b-1-d\geq 0;\ d \leq b-1$.
Altogether, $3 \leq d \leq b-1$.  Also, by our basic principle
above one can easily verify that indeed $(d_{2,b})^2 + (d^2)^2 = $
integer square.  This then establishes the first family.

\vspace{.15in}

\noindent {\bf Family F$_1$}  \fbox{\parbox{4.0in}{The following
family of Type 2 triangles is described in terms of two
independent positive integer parameters $\ell$ and $q$ which
satisfy the condition $\ell^2 \leq 2q^2 -2$.  This family consists
of all Type 2 triangles of the form $T_2(2,b,d)$, where the
integers $b$ and $d$ are defined as follows:  $b = 2rq-1,\
d=r^2-q^2,\ r = q+\ell$.}}

\vspace{.15in}

\noindent To construct the second family, let $\ell$ and $q$ be
positive integers such that, this time, $\ell^2 \geq 2q^2+2$.
Also, let $r=q+\ell,\ b=r^2-q^2 -1,\ d = 2rq$.  By inspection, $d
\geq 2$.  Clearly, $\ell^2 \geq 2q^2 + 2 > 2q^2$, which shows that
$\ell > q$.  So, $b=r^2 - q^2 -1 = (\ell + q)^2 - q^2 - 1 \geq
2q\ell \geq 4$, since $\ell > q$.  Thus $b \geq 4$.  Moreover,
$b-1-d=r^2-q^2 -1-1-2rq = (q+\ell)^2 - q^2 -2 -2 (q + \ell)q =
\ell^2 - 2 -2q^2 \geq 0$.  Thus, altogether $2 \leq d \leq b-1$;
and by our basic principle above, we also have
$(d_{2,b})^2+(d^2)^2 = $ integer square.

\vspace{.15in}

\noindent {\bf Family F$_2$} \fbox{\parbox{4.0in}{The following
family of Type 2 triangles is described in terms of two
independent positive integer parameters which satisfy the
condition $\ell^2 \geq 2q^2 + 2$.  This family consists of all
Type 2 triangles of the form $T_2(2,b,d)$, where the integers $b$
and $d$ are defined as follows: \linebreak $b=r^2-q^2-1,\ d=2rq,\
r = q+\ell$. }}

\vspace{.15in}

\item[B.] {\it Three families of Type 1 Triangles}

Let $r,q$ be positive integers with $r > q$.  Let $b=r^2+q^2 -1,\
d=r^2-q^2$.  Then $d \geq (q+1)^2-q^2 = 2q+1 \geq 3$; and let
$b+1-d = 2q^2 \geq 2$, so that $3 \leq d \leq b-1$.  Also, by our
basic principle, $(b+1)^2 - d^2 = (2rq)^2$ and thus $(d_{2,b})^2 -
(d^2)^2 = $ square.

\vspace{.15in}

\noindent {\bf Family S$_1$} \fbox{\parbox{4.0in}{The following
family of Type 1 triangles is described in terms of two
independent positive integer parameters $r$ and $q$ satisfying $r
> q$.  This family consists of all Type 1 triangles of the form
$T_1(2,b,d)$ where $b$ and $d$ are defined as follows: \linebreak
$b = r^2+q^2-1,\ d=r^2-q^2,\ r
> q \geq 1$.}}

\vspace{.15in}

Next, let $r,q$ be positive integers with $r \geq q + 2$.  Take
$b=r^2+q^2-1$, \linebreak $ d = 2rq$.  Then, clearly $d \geq 2$
and $b+1 -d=(r-q)^2 \geq 4$, so that $2 \leq d \leq b -3 < b-1$.
As before, by our basic principle, we have $(d_{2,b})^2 - (d^2)^2
= $ square.

\vspace{.15in}

\noindent {\bf Family S$_2$} \fbox{\parbox{4.0in}{The following
family of Type 1 triangles is described in terms of two
independent positive integer parameters.  This family consists of
all Type 1 triangles of the form $T_1(2,b,d)$; where $b$ and $d$
are defined as follows:   \linebreak $b=r^2+q^2 -1,\ d=2rq,\ r
\geq q+2$.}}

\vspace{.15in}

Finally, let $b \geq 4$, with $b$ being an integer square,
$b=k^2,\ k \in {\mathbb Z}^+$.  Let $d = b-1$.  Then $(d_{2,b})^2
- (d^2)^2 = $ square, since $(b+1)^2 - d^2 = (b+1)^2 - (b-1)^2 =
4b = 4k^2$.

\vspace{.15in}

\noindent {\bf Family U}  \fbox{\parbox{4.0in}{The following
family of Type 1 triangles is described in terms of one positive
integer parameter $t$.  This family consists of all Type 1
triangles of the form $T_1(2,b,d)$ where $b = t^2,\ d=b-1=t^2-1, \
t \geq 2, \ t \in {\mathbb Z}^+$.}}
\end{enumerate}

\section{Corollaries}

\noindent {\it Corollaries of Family $S_1$}

Let $d$ be an odd integer with $d \geq 3$. Then, there exists a
triangle which is a $T_1(2,b,d)$ for some positive integer $b$
with $d\leq b-1$.

\vspace{.15in}

\noindent {\bf Proof:}

We have $d = 2v+1$, for some integer $v \geq 1$.  Set $u = v+1$,
so that $d=u^2 -v^2$.  By Family S$_1$, if we take $b = u^2 + v^2
-1 = \dfrac{d^2-1}{2}$ then $d \leq b-1$ and $(d_{2,b})^2 -
(d^2)^2 = $ square \hfill $\Box$

\vspace{.15in}

\noindent {\it Corollaries of Family S$_2$}

Let $d$ be an even integer with $d \geq 6$.  Then, there exists a
triangle which is $T_1(2,b,d)$ for some positive integer $b$ with
$d \leq b-1$.

\noindent {\bf Proof:}

We have $d = 2v$, for some integer $v \geq 3$.  Let $u =1 $, so
that $v-u \geq 2$.  Then, if we take $b=v^2 = \dfrac{d^2}{4}$, by
Family S$_2$, we have $d \leq b-1$ and $(d_{2,b})^2 - (d^2)^2 = $
square. \hfill $\Box$

\vspace{.15in}

\noindent {\it Remark:}  It is easy to show that for $d=2,4$ there
exists no integer $b$ with $d \leq b -1$ and $(d_{2,b})^2 -
(d^2)^2 =$ square.

\vspace{.15in}

\noindent {\it Corollary of Family F$_1$}

Let $d$ be an odd integer with $d \geq 5$.  Then, there exists a
triangle which is $T_2(2,b,d)$ for some positive integer $b$ with
$d \leq b-1$.

\vspace{.15in}

\noindent {\bf Proof:}

We have $d = 2q+1$, for some positive integer $q \geq 2$.

Let $\ell = 1$, so that $\ell^2 \leq 2q^2 -2$.  Also, note that if
we take $r=q+\ell,\ r = q+1$, then $ d=r^2-q^2 = (q+1)^2 - q^2$.
Therefore, by Family F$_1$, if we set $b=2(q+1)q-1 = 2rq-1 =
\dfrac{d^2-3}{2}$.  It follows that $d \leq b-1$ and $(d_{2,b})^2+
(d^2)^2 = $ square. \hfill $\Box$

\vspace{.15in}

\noindent {\it Corollary to Family F$_2$}

Let $d$ be an even integer with $d \geq 6$.  Then, there exists a
triangle which is $T_2(2,b,d)$ for some positive integer $b$ with
$d \leq b-1$.

\vspace{.15in}

\noindent {\bf Proof:}  We have $d = 2r$, for some integer $r \geq
3$.  Let $q=1$, so that $d=2rq$; and $\ell = r-q = r-1 \geq 2$.
Then $\ell^2 \geq 2q^2 + 2 = 4$.  From Family F$_2$, it is now
clear that if we set $b = r^2-2 = \dfrac{d^2-8}{4}$, then $d \leq
b-1$ and $(d_{2,b})^2 + (d^2)^2 = $ square.  \hfill $\Box$

Note that the corollaries of F$_1$ and F$_2$ complement the result
of Theorem 5.

\end{document}